\documentclass[11pt,a4paper,twoside,reqno,nosumlimits]{amsart}

\usepackage{graphicx}
\usepackage{epstopdf}
\usepackage{amsmath}
\usepackage{amsfonts}
\usepackage{bbm}
\usepackage{braket}

\usepackage[normalem]{ulem}

\usepackage{color}
\usepackage{subfigure}
\usepackage[margin=2.0cm]{geometry}
\makeatletter

\newcommand{\Rmnum}[1]{\expandafter\@slowromancap\romannumeral #1@}

\makeatother

\begin{document}
\title{Exploring simplicity bias in 1D dynamical systems}

\author{Kamal Dingle$^{1,2}$, Mohammad Alaskandarani$^{1}$, Boumediene Hamzi$^{3,4}$, Ard A. Louis$^2$}

\address{$^1$Centre for Applied Mathematics and Bioinformatics,
Department of Mathematics and Natural Sciences, 
Gulf University for Science and Technology, Kuwait}
\email{dingle.k@gust.edu.kw}

\address{$^2$Rudolf Peierls Centre for Theoretical Physics$\text{,}$\\ University of Oxford$\text{,}$ Parks Road,
Oxford, OX1 3PU, UK}
\email{ard.louis@physics.ox.ac.uk}

 \address{$^3$ Department of Computing and Mathematical Sciences, Caltech, CA, USA.}
\address{$^4$ The Alan Turing Institute, London, UK.}
\email{boumediene.hamzi@gmail.com}, 

\date{\today}

\begin{abstract}
\noindent 
Arguments inspired by algorithmic information theory predict an inverse relation between the probability and complexity of output patterns in a wide range of input-output maps. This phenomenon is known as \emph{simplicity bias}. By viewing the parameters of dynamical systems as inputs, and resulting (digitised) trajectories as outputs, we study simplicity bias in the  logistic map, Gauss map, sine map, Bernoulli map, and tent map. 
We find that the logistic map, Gauss map, and sine map all exhibit simplicity bias upon sampling of map initial values and parameter values, but the Bernoulli map and tent map do not. The simplicity bias upper bound on output pattern probability is used to make \emph{a priori} predictions for the probability of output patterns. In some cases, the predictions are surprisingly accurate, given that almost no details of the underlying dynamical systems are assumed.  More generally, we argue that studying probability-complexity relationships  may be a useful tool in studying patterns in dynamical systems.
\\

\end{abstract}

\maketitle

\section{Introduction}


In recent years several studies of \emph{simplicity bias} have been made in input-output maps, in which a general inverse relationship between the complexity of outputs and their respective probabilities has been observed \cite{dingle2018input,dingle2020generic}. More specifically, using arguments inspired by   \emph{algorithmic information theory} \cite{solomonoff1960preliminary,kolmogorov1965three,chaitin1975theory} (AIT), and specifically \emph{algorithmic probability}, an upper bound on the probability $P(x)$ of observing output pattern $x$ was presented \cite{dingle2018input}, with the bound depending on the estimated Kolmogorov complexity of $x$. The upper bound implies that complex output patterns must have low probabilities, while high probability outputs must be simple. Example systems where simplicity bias has been observed include RNA structures \cite{dingle2018input,dingle2023multiclass}, differential equation solutions \cite{dingle2018input}, finite state transducers \cite{dingle2020generic}, time series patterns in natural data \cite{dingle2023note}, natural protein structures \cite{johnston2022symmetry},  among others.

A full understanding of exactly which systems will, and will not, show simplicity bias is still lacking, but the phenomenon is expected to appear in a wide class of input-output maps, under fairly general conditions. Some of these conditions were suggested in ref.\ \cite{dingle2018input}, including (1) that the number of inputs should be much larger than the number of outputs, (2) the number of outputs should be large, and (3) that the map should be `simple' (technically of $O(1)$ complexity) to prevent the map itself from dominating over inputs in defining output patterns. Indeed, if an arbitrarily complex map was permitted, outputs could have arbitrary complexities and probabilities, and thereby remove any connection between probability and complexity.  Finally (4), because many AIT applications rely on approximations of Kolmogorov complexity via standard lossless compression algorithms \cite{lempel1976complexity,ziv1977universal} (but see \cite{delahaye2012numerical,soler2014calculating} for a fundamentally different approach), another condition proposed is that the map should not generate pseudo-random outputs like $\pi=3.1415...$ which standard compressors cannot handle effectively. The presence of such outputs may yield high probability outputs which appear `complex' hence apparently violating simplicity bias, but which are in fact simple. 

To explore further the presence of simplicity bias in dynamical systems and physics, and also in maps which test the boundaries of the conditions for simplicity bias described above,   we examine the output trajectories of a selection of 1D maps from chaos theory, namely the logistic map, the Gauss (``mouse'') map, the sine map, the Bernoulli map, and the tent map. Starting with its popularisation by Robert May \cite{may1976simple} in the 1970's, the logistic map has been heavily studied.  The map is a textbook example of a simple system which can produce simple, complex, and chaotic patterns via iterations of the map \cite{hasselblatt2003first}. This map is related to a very wide range of nonlinear models with applications in epidemiology, economics, physics, time series, etc. \cite{hilborn2000chaos}. Due to the popularity of the logistic map, and because its trajectory outputs can depict simple as well as complex, chaotic, and even pseudo-random patterns, we  focus primarily on the logistic map, but we also analyse the other mentioned maps. 

Although not restricted to studying binary strings, most AIT results are framed in the context of binary strings and hence applications are easier in the same framework. Thus in this work we will study simplicity bias in digitised binary trajectories of these example 1D dynamical systems.  Our main conclusions are that simplicity bias appears in the logistic map, the Gauss (``mouse'') map, and the sine map, and hence we suggest that simplicity bias may also appear in natural dynamical systems more generally. As a broader context of motivation, this work contributes to research is the intersection of dynamical systems, and AIT, and machine learning.

\section{Background and problem set up}

\subsection{Background theory and pertinent results}
\subsubsection{AIT and Kolmogorov complexity}
We give some basic background regarding AIT now, and describe simplicity bias in more detail. Note that the current work will not involve detailed AIT or related theory, so we only give a brief survey of relevant results without giving many formal details here. There are many standard texts which the interested reader can refer to if needed, e.g., refs.\ \cite{li2008introduction,calude2002information,gacs1988lecture,shen2022kolmogorov}.

Within computer science, \emph{algorithmic information theory} \cite{solomonoff1960preliminary,kolmogorov1965three,chaitin1975theory} (AIT) directly connects computation, computability theory, and information theory. The central quantity of AIT is \emph{Kolmogorov complexity}, $K(x)$, which measures the complexity of an individual object $x$ as the amount of information required to describe or generate $x$. $K(x)$ is more technically defined as the length of a shortest program which runs on an (optimal prefix) \emph{universal Turing machine} (UTM) \cite{turing1936computable}, generates $x$, and halts. 
More formally, the Kolmogorov complexity $K_U(x)$ of a string $x$ with respect to $U$,  is defined \cite{solomonoff1960preliminary,kolmogorov1965three,chaitin1975theory} as
\begin{equation}
K_U(x) = \min_{p}\{|p|: U(p)=x\}
\end{equation}
where $p$ is a binary program for a prefix optimal UTM $U$, and $|p|$ indicates the length of the program $p$ in bits.  Due to the invariance theorem \cite{li2008introduction} for any two optimal UTMs $U$ and $V$, $K_U(x) = K_V(x)+O(1)$ so that the complexity of $x$ is independent of the machine, up to additive constants. Hence  we conventionally drop the subscript $U$ in $K_U(x)$, and speak of `the' Kolmogorov complexity $K(x)$. Informally, $K(x)$ can be defined as the length of a shortest program that produces $x$, or simply as the size in bits of the compressed version of $x$ (assuming a perfect compressor). If $x$ contains repeating patterns like $x=1010101010101010$ then it is easy to compress, and hence $K(x)$ will be small. On the other hand, a randomly generated bit string of length $n$ is highly unlikely to contain any significant patterns, and hence can only be described via specifying each bit separately without any compression, so that $K(x)\approx n$ bits. Other more expressive names for $K(x)$ are \emph{descriptional complexity}, \emph{algorithmic complexity}, and \emph{program-size complexity}, each of which highlight the idea that $K(x)$ is measuring the amount of information to describe or generate $x$ precisely and unambiguously. Note that Shannon information and Kolmogorov complexity are related \cite{grunwald2004shannon}, but differ fundamentally in that Shannon information quantifies the information or complexity of a random source, while Kolmogorov complexity quantifies the information of individual sequences or objects.

An increasing number of studies show that AIT and Kolmogorov complexity can be successfully applied in physics, including thermodynamics \cite{bennett1982thermodynamics,kolchinsky2020thermodynamic,zurek1989algorithmic,kolchinsky2023generalized}, quantum physics \cite{mueller2020law}, and entropy estimation \cite{avinery2019universal,martiniani2019quantifying}. Further, applications to biology \cite{ferragina2007compression,johnston2022symmetry,adams2017formal}, other natural sciences \cite{devine2020algorithmic}, and engineering, are also numerous \cite{vitanyi2013similarity,li2008introduction,cilibrasi2005clustering}. 
 
 \subsubsection{The coding theorem and algorithmic probability}
 
An important result in AIT is Levin's \cite{levin1974laws} coding theorem, establishing a fundamental connection between $K(x)$ and probability predictions. Mathematically, it states that 
\begin{equation}
P(x) = 2^{-K(x)+O(1)}\label{eq:CD}
\end{equation}
where $P(x)$ is the probability that an output $x$ is generated by a (prefix optimal) UTM fed with a random binary program. Thus, high complexity outputs have exponentially low probability, and simple outputs must have high probability. This is a profound result which links notions of data compression and probability in a direct way. $P(x)$ is also known as the \emph{algorithmic probability} of $x$. Given the broad reaching and striking nature of this theorem, it is somewhat surprising that it is not more widely studied in the natural sciences. The reason in part for this inattention is that AIT results are often difficult to apply directly in real-world contexts, due to a number of issues including the fact that $K(x)$ is formally uncomputable and the ubiquitous use of UTMs.

 \subsubsection{The simplicity bias bound}

Coding theorem-like behaviour in real-world input-output maps has been studied recently, leading to the observation of a phenomenon called \emph{simplicity bias} \cite{dingle2018input} (see also ref.\ \cite{buchanan2018natural}). Simplicity bias is captured mathematically as
\begin{equation}
P(x)\leq 2^{-a\tilde{K}(x)-b}\label{eq:simplicity bias}
\end{equation}
where $P(x)$ is the (computable) probability of observing output $x$ on random choice of inputs, and $\tilde{K}(x)$ is the approximate Kolmogorov complexity of the output $x$: complex outputs from input-output maps have lower probabilities, and high probability outputs are simpler. The constants $a>0$ and $b$ can be fit with little sampling and often even predicted without recourse to sampling \cite{dingle2018input}. We will assume that $b=0$ in Eq.\ (\ref{eq:simplicity bias}) throughout this work, which is a default assumption as argued and discussed in ref.\ \cite{dingle2018input}. There is also a conditional version of the simplicity bias equation \cite{dingle2022predicting}.

The ways in which simplicity bias differs from Levin's coding theorem include that it does not assume UTMs, uses approximations of complexities, and for many outputs $P(x)\ll 2^{-K(x)}$. Hence the abundance of low complexity, low probability outputs \cite{dingle2020generic,alaskandarani2023low} is a signature of simplicity bias.

 \subsubsection{Estimating pattern complexity} 

To estimate complexity, we follow ref.\ \cite{dingle2018input}  and use
\begin{equation}
C_{LZ}(x) =\begin{cases}
     \log_2(n), &  \hspace*{-0.3cm}  \text{$x=0^n$ or $1^n$}\\
    \log_2(n) [N_w(x_1...x_n) + N_{w}(x_n...x_1)]/2, & \hspace*{-0.2cm} \text{otherwise}
  \end{cases}\label{eq:CLZ}
\end{equation}
where $N_w(x)$ comes from the 1976 Lempel and Ziv complexity measure  \cite{lempel1976complexity}, and 
where the simplest strings $0^n$ and $1^n$ are separated out because  $ N_{w}(x)$ assigns complexity $K=1$ to the string 0 or 1, but complexity 2 to $0^n$ or $1^n$ for $n\geq2$, whereas the true Kolmogorov complexity of such a trivial string actually scales as $\log_2(n)$ for typical $n$, because one only needs to encode $n$. 
Having said that, the minimum possible value is $K(x)\approx0$ for a simple set, and so e.g. for binary strings of length $n$ we can expect $0 \leq K(x) \leq n$ bits. Because for a random string of length $n$ the value $C_{LZ}(x)$ is often much larger than $n$, especially for short strings, we scale the complexity so that $a$ in Eq.\ (\ref{eq:simplicity bias}) is set to $a=1$ via
\begin{equation}
\tilde{K}(x) = \log_2(M) \cdot \frac{ C_{LZ}(x) - \min_x (C_{LZ})}{\max_x (C_{LZ}) - \min_x (C_{LZ}) } \label{eq:Kscaled}
\end{equation}
where $M$ is the maximum possible number of output patterns in the system, and the min and max complexities are over all strings $x$ which the map can generate.  $\tilde{K}(x)$ is the approximation to Kolmogorov complexity that we use throughout. 
This scaling results in $0\leq \tilde{K}(x) \leq n$ which is the desirable range of values.

\begin{figure}[htp]
\begin{center}
\subfigure[]{\label{fig:edge-a}\includegraphics[height=7cm,width=9cm]{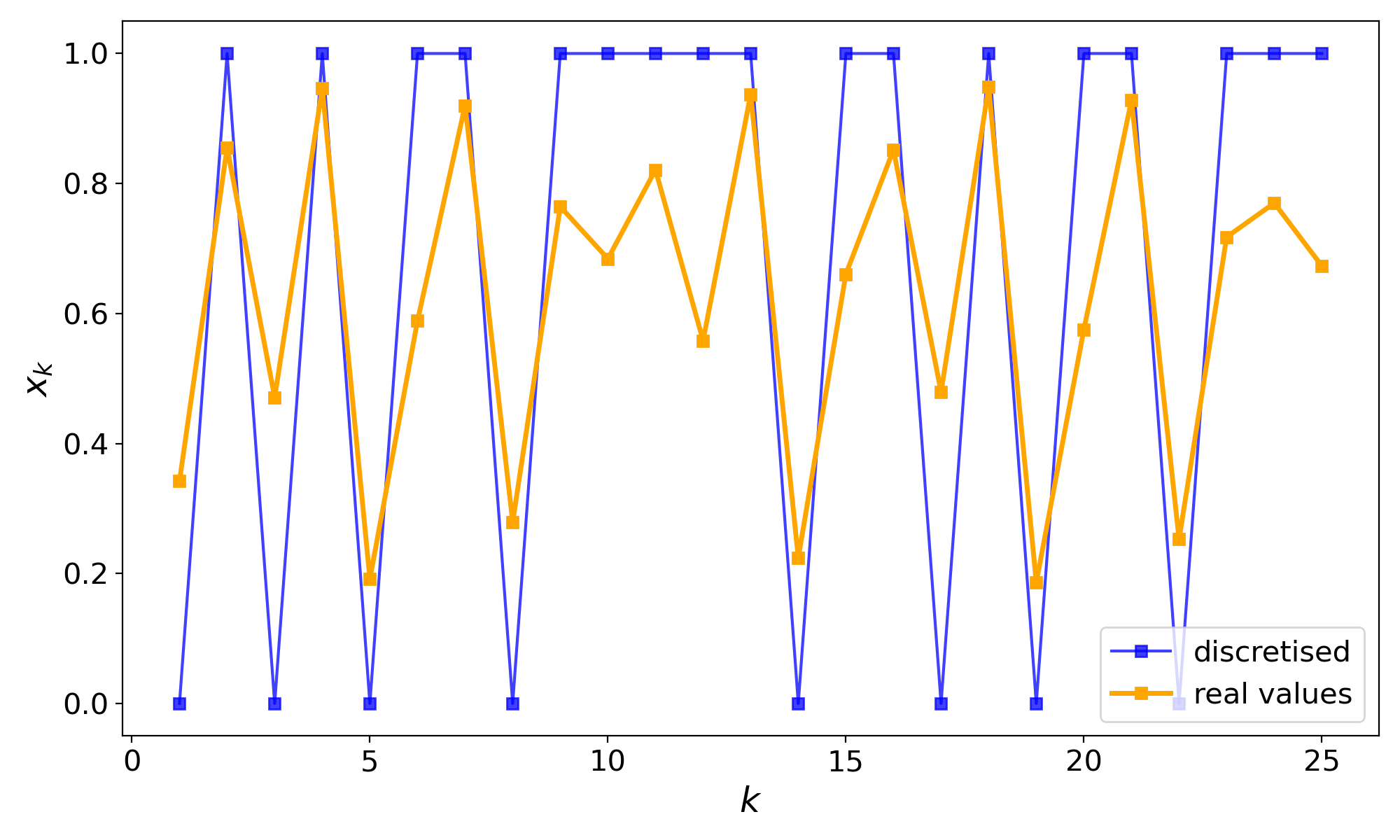}}
\end{center}
\caption{An example of a real-valued (orange) and digitised (blue) trajectory of the logistic map, with $\mu=3.8$ and $x_0=0.1$. The discretisation is defined by writing 1 if $x_k\geq 0.5$ and 0 otherwise, resulting in the pattern $x=$ 0101011011111011010110111 which has length $n=$ 25 bits. }
\label{fig:discretisation}
\end{figure}

\subsection{Digitised map trajectories}
The AIT coding theorem is framed in terms of random inputs or `programs' and resultant output strings or patterns. The core idea of the simplicity bias bound is that for a large range of systems, uniformly sampling input parameter values (`programs') yields an exponential variation in output pattern probabilities, with high probability patterns having low complexities. While dynamical systems may not appear to fit this input-output framework, they can be viewed as input-output functions if the initial values and parameters are considered input `programs' which are then computed to generate `output' dynamical system trajectories. Because output pattern complexities and probabilities are more easily calculated if outputs are binary strings (and as above, AIT is mainly framed in terms of binary strings), we will digitise the real-valued trajectories into binary strings.

To illustrate this input-output framework, consider binary sequence trajectories resulting from digitised realisations of the logistic map
\begin{equation}
x_{k+1}=\mu x_k(1-x_k) \label{eq:logisticmap}
\end{equation}
where the inputs are the pair of values $x_0\in$(0.0, 1.0) and $\mu\in$ (0.0, 4.0]. For the outputs, the map is first iterated to obtain a sequence of $n$ real values $x_1,x_2,x_3,\dots,x_n$ in [0,1]. Similar to the field of symbolic dynamics \cite{lind1995introduction}, the real-valued trajectory is digitised to become a binary string output sequence $x$ by applying a threshold, and writing 1 if $x_k\geq0.5$ and 0 otherwise \cite{kanso2009logistic}. Hence a binary sequence output $x$ of $n$ bits is generated for each input pair $\mu$ and $x_0$. 

By way of example, consider choosing $\mu=3.8$ and $x_0=0.1$ with $n=25$, then from iterating Eq.\ (\ref{eq:logisticmap}) with $k=1,2,...,25$ we obtain the real-valued trajectory 
\[x_1, x_2,\dots, x_{25} = 0.34, 0.86, \dots , 0.25, 0.72, 0.77, 0.67\]
which after digitisation becomes the binary string
\[x=0101011011111011010110111\]
 Figure \ref{fig:discretisation}  illustrates the trajectory and digitisation procedure. The choice of $n$ is a balance between being long enough that many different patterns can emerge, and short enough that decent frequency and probability estimates can be made without the need for excessive sampling. Different $\mu$ and $x_0$ inputs can yield different binary strings $x\in\{0,1\}^n$.

Note that throughout the paper, we will ignore the first 50 iterations of a map, so as to discard the transient dynamics of the iterated maps. This step of ignoring the initial iterations is common practice \cite{berger2001chaos}.

\begin{figure*}[htp]
\begin{center}
\subfigure[]{\label{fig:edge-a}\includegraphics[height=7cm,width=8cm]{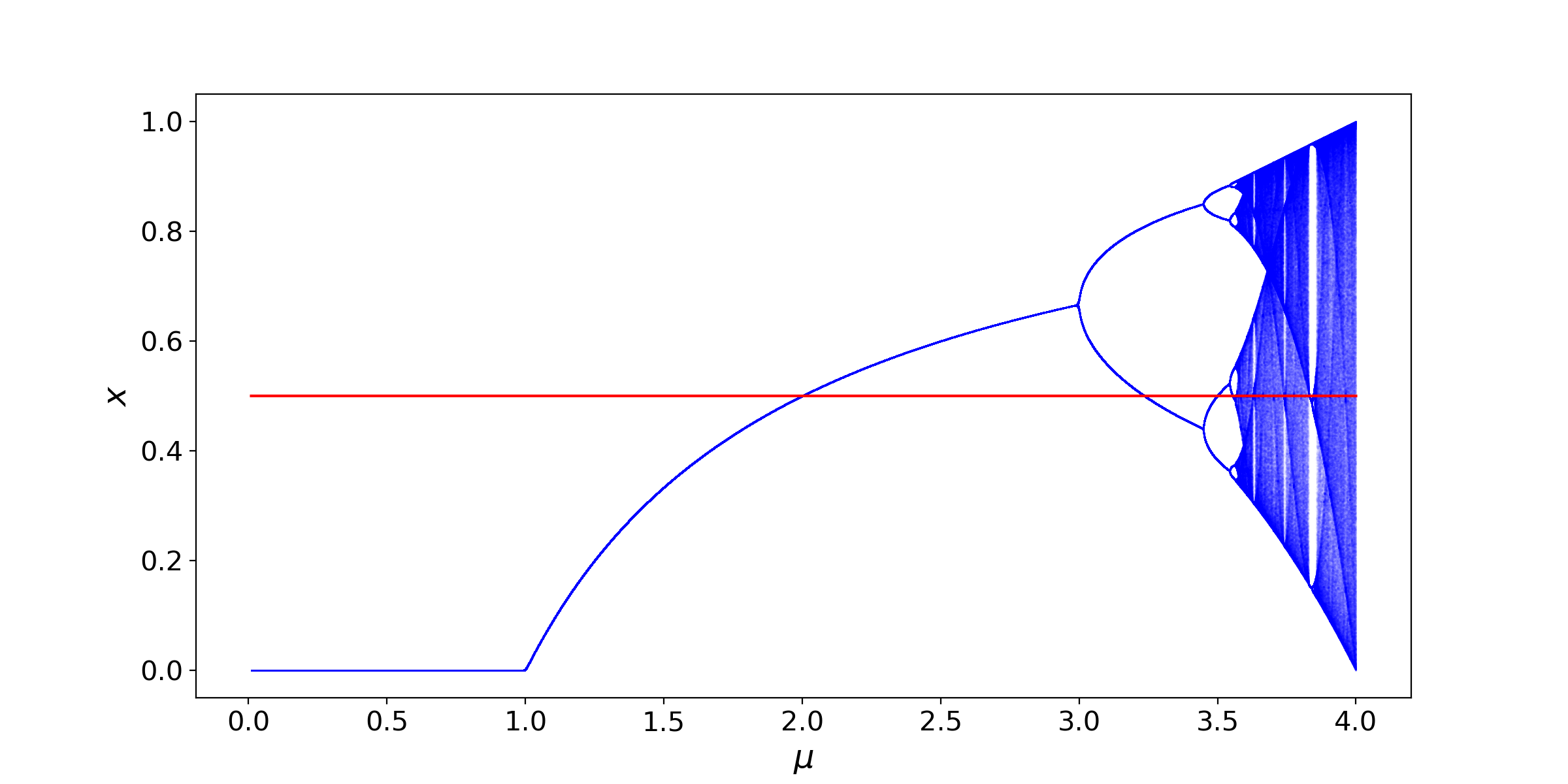}}
\subfigure[]{\label{fig:edge-a}\includegraphics[height=7cm,width=8cm]{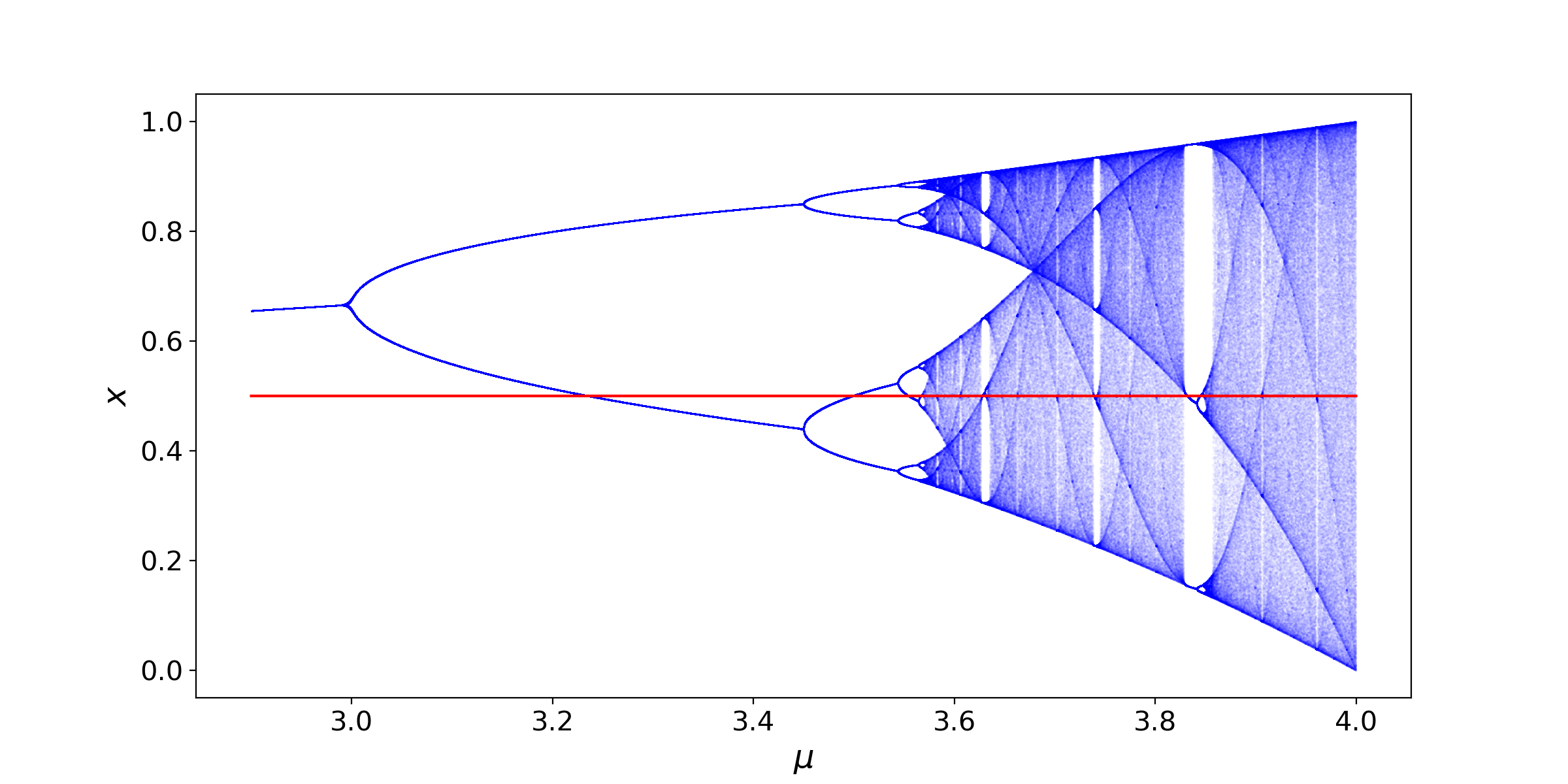}}
\end{center}
\caption{Bifurcation diagram for the logistic map. In (a), the diagram for parameters $\mu\in$(0, 4.0]; and in (b), for values $\mu\in$(2.9, 4.0]. The value 0.5 has been highlighted in red, to indicate the cut off threshold used to digitise trajectories by a value of $0$ if the output is below the threshold, and a value of $1$ if it is greater than or equal to the threshold.}
\label{fig:bifurcation}
\end{figure*}

\section{Results}

\subsection{Logistic map}
The logistic map in Eq.\ (\ref{eq:logisticmap}) is possibly the best known and most iconic 1D map studied in dynamical systems  and chaos theory.  To our knowledge, there have not been many studies of  the complexity of digitised trajectories of the logistic map. A notable exception is the work by  Kaspar and Schuster \cite{kaspar1987easily} who also studied the estimated complexity of resulting binary string trajectories. However their work was fundamentally different in that it was not concerned with simplicity bias or in  estimating the probability of different outputs.

\subsubsection{Parameter intervals}

The intervals  $x_0\in$ (0.0, 1.0) and $\mu\in(0.0, 4.0]$ are standard ranges in which the logistic map is studied and it can be shown \cite{berger2001chaos}, for example, that for $\mu>4.0$ almost all trajectories are not confined to [0,1]; similarly if $x_0\notin(0,1)$ the behaviour can be unbounded or trivial.  For some large values of $\mu$ (e.g., $\mu=4.0$), almost all initial values $x_0$ yield complex and chaotic outputs \cite{hasselblatt2003first}, and the distribution over digitised trajectories is roughly uniform, with little bias. Note that we use the word `bias' to describe a strongly non-uniform distribution of the probability to obtain a given binary string output. 

In Figure \ref{fig:bifurcation}(a) the bifurcation diagram of the logistic map is shown, in which for different values of $\mu$ the asymptotic $x_k$ values are depicted. The diagram shows fixed points, oscillations, and non-period behaviour.
Also added is the value $0.5$ as a red line, highlighting the digitising threshold. It is known  \cite{hasselblatt2003first} that if $\mu\in[0.0,1.0]$, then the trajectories tend to 0. Because we truncate at 0.5, therefore the corresponding binary string would be $x=...0000$. If $\mu\in(1.0,2.0]$ then the trajectories tend to $1-(1/\mu)$, which means that we expect to see binary strings $x=...0000$. For $\mu\in(2.0,3.0]$, the resulting pattern would be $x=...1111$ because $1-(1/\mu)$ is greater than 0.5. For $\mu$ from 3.0 to  $\approx$3.3, we still expect $x=...1111$ because although the first bifurcation appears at $\mu=$3.0, and oscillations with period two begin, until about 3.3 both values of the oscillations are larger than 0.5. Figure \ref{fig:bifurcation}(b) shows the same bifurcation diagram, but zoomed in on the larger values of $\mu$. From about 3.3 to about 3.5, oscillations between two or four values appear, and larger values of $\mu$ can yield oscillation periods of 8, 16, etc., which will yield patterns such as $x=0101...$. 
Chaotic trajectories do not occur until \cite{hasselblatt2003first} $\mu\geq3.56994567...\approx 3.57$ (\texttt{https://oeis.org/A098587}), and so  the $\mu$ interval 3.0 to about 3.5699 contains  period-doubling bifurcations, in which oscillations of exponentially increasing frequency occur but no truly chaotic trajectories. Finally, for $\mu$ between  $\approx$3.57 and 4.0, more complex and chaotic patterns can emerge, but  also `islands' of stability appear in this interval, as can be observed in Figure \ref{fig:bifurcation}(b). 

In our numerical simulations, we will separately investigate various intervals for $\mu$, namely $\mu$ sampled uniformly from (0.0, 4.0], [3.0, 4.0], [3.57, 4.0], and finally also fixing $\mu=4.0$. The motivation for choosing these intervals is as follows: Given that the interval  $\mu\in(0.0, 4.0]$ is the standard range in which the map is studied, the most `natural' sampling strategy is to sample uniformly across this interval. As for choosing intervals [3.0, 4.0] and [3.57, 4.0], these are interesting test cases because we can expect to see complex patterns appearing more frequently in these intervals. Finally, when fixing $\mu=4.0$ most trajectories will be highly complex, and so we might expect simplicity bias to disappear. 
\subsubsection{Connection of simplicity and probability}

Reflecting on the bifurcation diagram and above related comments,  it is easily seen that by uniformly sampling $\mu$ from (0.0, 4.0], some simple binary patterns will have high probability, and highly complex patterns will have low probability just because most of the interval (0.0, 4.0] yields fixed points or low period oscillations. Hence, some general form of bias towards simplicity is expected for the logistic map. However, what is not \emph{a priori} obvious is whether or to what extend the simplicity bias bound in Eq.\ (\ref{eq:simplicity bias}) will be followed or have predictive value.

\begin{figure*}[htp]
\begin{center}
\subfigure[]{\label{fig:edge-a}\includegraphics[height=5.5cm,width=5.5cm]{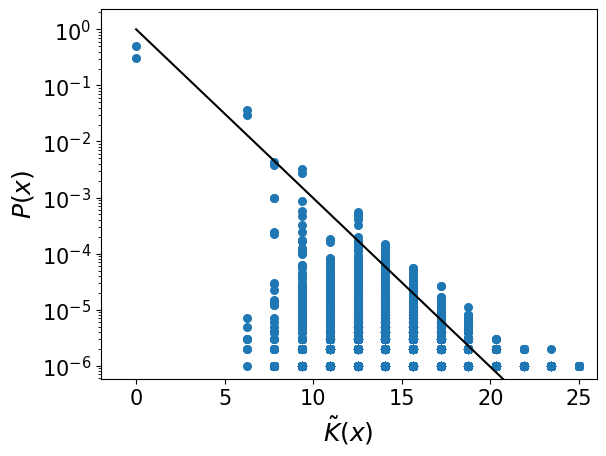}}
\subfigure[]{\label{fig:edge-a}\includegraphics[height=5.5cm,width=5.5cm]{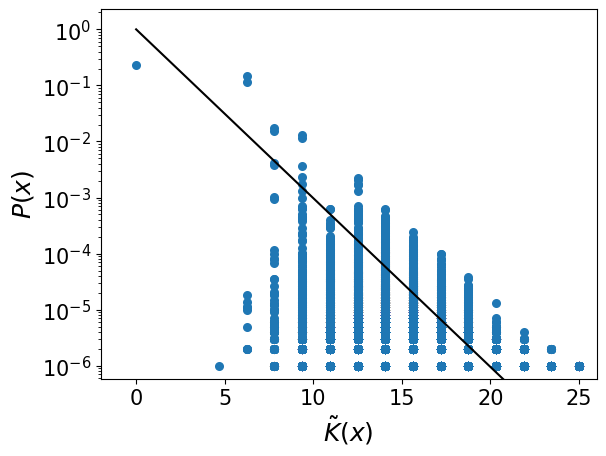}}\\
\subfigure[]{\label{fig:edge-a}\includegraphics[height=5.5cm,width=5.5cm]{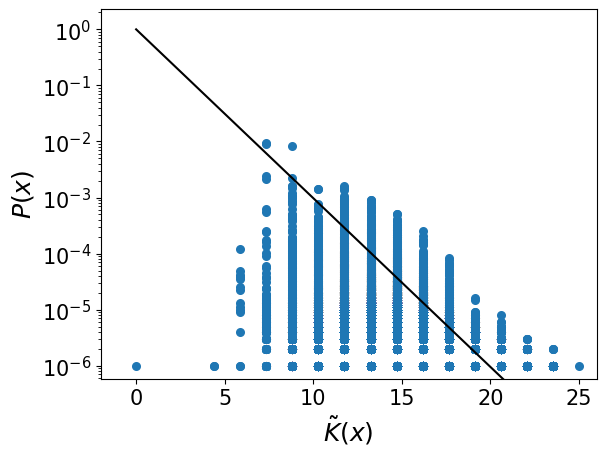}}
\subfigure[]{\label{fig:edge-a}\includegraphics[height=5.5cm,width=5.5cm]{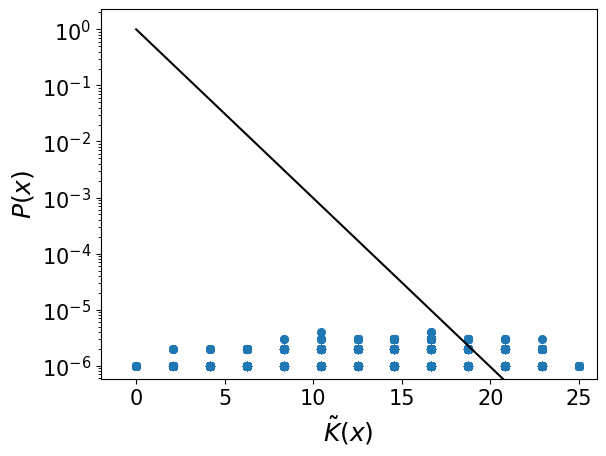}}
\end{center}
\caption{Simplicity bias in the digitised logistic map from random samples with $x_0\in(0,1)$ and $\mu$ sampled in different intervals. Each blue datapoint corresponds to a different binary digitised trajectory $x$ of length 25 bits. The black line is the upper bound prediction of Eq.\ (\ref{eq:simplicity bias}). (a) Clear simplicity bias for $\mu\in$(0.0, 4.0] with $P(x)$ closely following the upper bound, except for low frequency and high complexity outputs which suffer from increased sampling noise; (b) simplicity bias is still present for  $\mu\in$[3.0, 4.0]; (c) the distribution of $P(x)$ becomes more flat (less biased) and simplicity bias is much less clear when $\mu\in$[3.57, 4.0] due to constraining the sampling to $\mu$-regions more likely to show chaos; (d) the distribution of $P(x)$ is roughly uniform when using $\mu=4.0$, with almost no bias, and hence no possibility of simplicity bias.   
}
\label{fig:bias_with_simplicity bias}
\end{figure*}

\subsubsection{Simplicity bias appears when bias appears}

Following the protocol for generating binary strings via logistic map trajectories described above, we now examine the probability $P(x)$ that the digitised trajectory $x$ is produced on random sampling of input parameters $\mu$ and $x_0$, done here with  $10^6$ samples over uniformly chosen random parameters.  Using Eq.\ (\ref{eq:simplicity bias}), we can make an \emph{a-priori} prediction for the upper bound decay upon sampling of the input parameters. (Recall that we ignore the first 50 iterations of the map, in order to exclude the transient dynamics.)

Figure \ref{fig:bias_with_simplicity bias}(a) shows that the upper bound prediction (black line) agrees remarkably well with the probability and complexity data for  different binary string output (blue dots);  we see that the logistic map displays simplicity bias when uniform sampling $\mu\in(0.0, 4.0]$, such that high probability outputs have low complexities, and complex outputs have low probability. The gradient prediction ($a=1$) of the black line used no information about the map, except that we assumed that almost all $2^n$ outputs of length $n$ bits are realisable. An upper-bound fit to the $\log_{10}P(x)$ data gives the slope as -0.18. Note that many output strings fall below the upper bound prediction as we expected from earlier studies \cite{dingle2018input,alaskandarani2023low}, but nonetheless it is known that randomly generated outputs tend to be close to the bound \cite{dingle2020generic}. Put differently, even though many output strings (blue dots) appear to be far below the bound, most of the probability mass for each complexity value is concentrated close to the bound. Thus, this simple bound  predicts $P(x)$ quite accurately by using the complexity values of output strings yet while otherwise completely ignoring the details of the underlying dynamical system.

In Figure \ref{fig:bias_with_simplicity bias}(b) we make a similar plot to panel (a) except that we restrict the sampling to $\mu\in[3.0, 4.0]$. The qualitative pattern in (b) is similar to (a), although simplicity bias is slightly less clear than in (a). An upper-bound fit to the $\log_{10}P(x)$ data gives the slope as -0.17. For Figure \ref{fig:bias_with_simplicity bias}(c)  the sampling interval was chosen so that $\mu\in [3.57, 4.0]$ which is the region containing some truly chaotic trajectories \cite{berger2001chaos}. Both bias and simplicity bias are less clear in this plot. An upper-bound fit to the $\log_{10}P(x)$ data gives the slope as -0.04. In Figure \ref{fig:bias_with_simplicity bias}(d) a plot is given in which $\mu$ is no longer sampled, but is fixed at $\mu=4.0$, so that almost all $x_0$ values yield chaotic trajectories. As expected, there is very little bias in the distribution, i.e., the distribution of $P(x)$ is roughly uniform, and hence no simplicity bias can be observed. An upper-bound fit to the $\log_{10}P(x)$ data gives the slope as 0.02. Figure \ref{fig:bias_with_simplicity bias}(d) still shows some very simple strings, which can be understood from the fact that $x_i\approx0$ $\Rightarrow x_{i+j}\approx 0$ for $j\approx 1$ even if the trajectories may become much more complex for $j\gg 1$. In other words, if the trajectory reaches a value close to zero, then it tends to stay close to zero for several iterations. These initial short trajectories with many 0s will have very low complexity. 

As argued above, the fact that simple patterns occur with high probability when sampling $\mu$ from (0.0, 4.0] is not in itself remarkable, and can be rationalised from the bifurcation diagram. 
However, what is noteworthy in our investigation is that we don't just  see a vague general inverse relation between complexity and probability (as the preceding discussion might predict), but rather that we see an exponential decay in probability with linearly increasing complexity which follows the upper bound of Eq.\ (\ref{eq:simplicity bias}) quite closely. 

See the Appendix \ref{numbiterations} for an illustration of the effects of using the different values of $n$.  See Appendix B for the same plots as in Figure 3, but with semi-transparent data points to highlight the distribution of data points.

\begin{figure*}[htp]
\begin{center}
\subfigure[]{\label{fig:edge-a}\includegraphics[height=5.5cm,width=5.5cm]{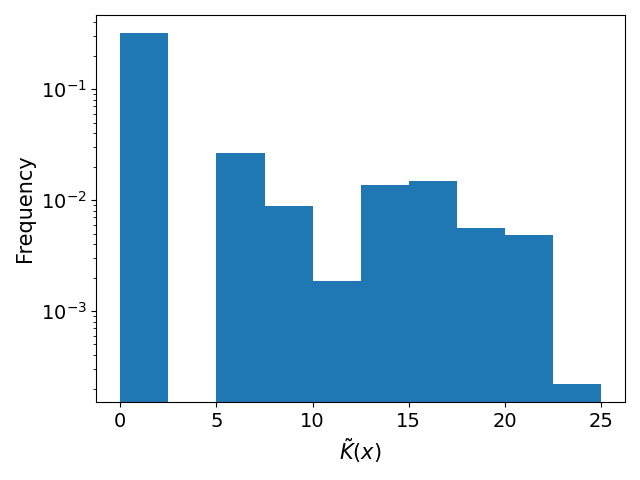}}
\subfigure[]{\label{fig:edge-a}\includegraphics[height=5.5cm,width=5.5cm]{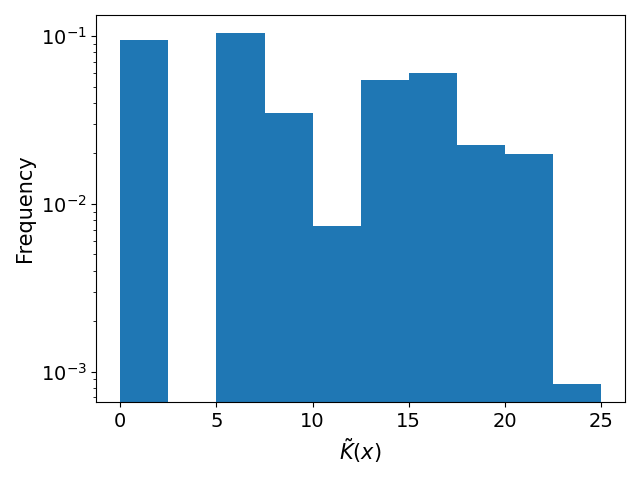}}\\
\subfigure[]{\label{fig:edge-a}\includegraphics[height=5.5cm,width=5.5cm]{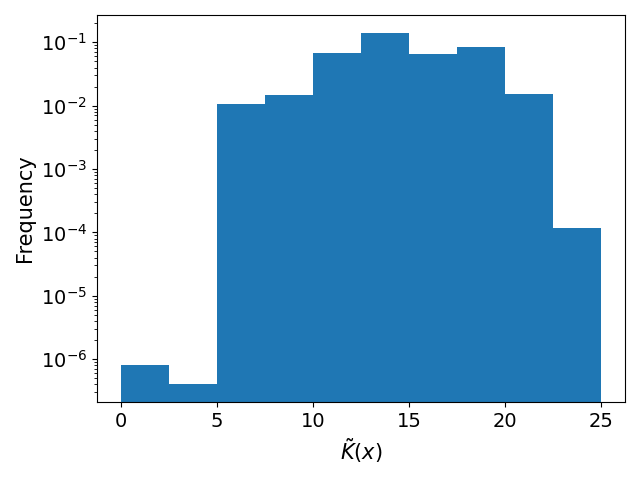}}
\subfigure[]{\label{fig:edge-a}\includegraphics[height=5.5cm,width=5.5cm]{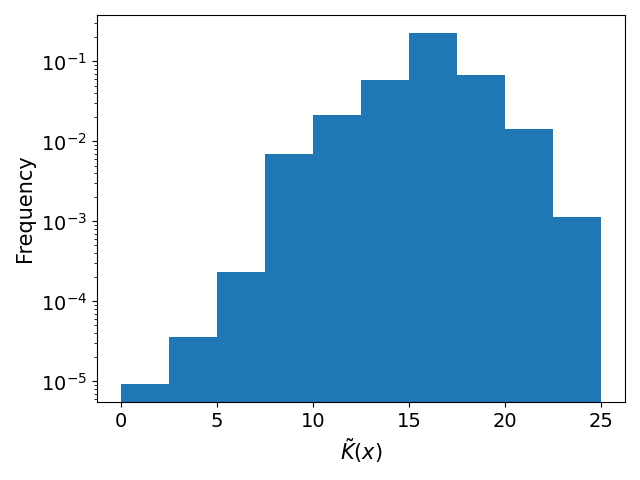}}
\subfigure[]{\label{fig:edge-a}\includegraphics[height=5.5cm,width=5.5cm]{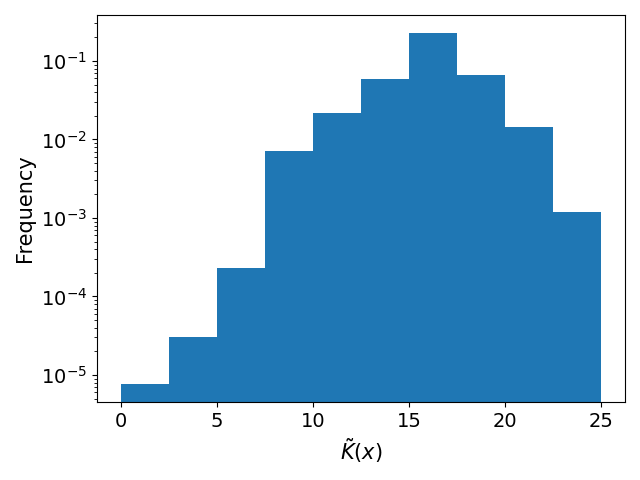}}
\end{center}
\caption{The distribution $P(\tilde{K}(x)=r)$  of output complexity values, with $x_0\in(0.0,1.0)$ and $\mu$ sampled from different intervals. (a) A roughly uniform complexity distribution for $\mu\in$(0.0, 4.0], with some bias towards lower complexities (mean is 3.4 bits); (b) Close to uniform distribution of complexities for $\mu\in$[3.0, 4.0], mean is 10.3 bits; (c) the distribution leans to higher complexities when $\mu\in$[3.57, 4.0], mean is 14.1 bits; (d) the distribution is  biased to higher complexities values when $\mu=4.0$ (mean is 16.4 bits); (e) for comparison, purely random binary strings of length 25 bits were generated (mean is 16.2 bits). The distributions of complexity values in (d) and (e) are very similar, but (a-c) show distinct differences. Calculating and comparing $P(K)$ is an efficient way of checking how simplicity biased a map is. 
}
\label{fig:distn_K}
\end{figure*}

\subsubsection{Distribution of complexities}

By itself, simplicity bias does not necessarily imply that the distribution of complexity values, $P(K(x)=r)$, will be skewed towards low complexities values. Rather, simplicity bias means that the individual probability of some simple strings will be higher than that of more complex strings. Because the probability of observing a given complexity $r$ depends on the product of the number of strings of complexity $r$ and their individual probabilities,  the distribution $P(K(x)=r)$ may in fact peak at high complexities, or low complexities, or have no peak at all possibly. 

To investigate this distribution in the logistic map, in Figure \ref{fig:distn_K} we plot the distribution of complexities for the different sampling intervals of $\mu$. As can be seen in Figure \ref{fig:distn_K}(a), when sampling from [0.0, 4.0] there is a bias towards lower complexities. As for (b), when $\mu$ is sampled from  [3.0, 4.0] the distribution of complexities is roughly uniform (at least on a log scale). In Figure \ref{fig:distn_K}(c), the distribution is also somewhat uniform but peaks around medium complexity values. In (d), there is a bias towards higher complexity values. 
For comparison, in (e) we also plot the distribution of complexities resulting from sampling purely random binary strings, and this distribution is very similar to that in (d) obtained when $\mu=4.0$.
It is noteworthy that in some of these cases, while there is some evidence of bias toward higher or lower complexities, the distributions still display spread and are not narrowly peaked (at least on a log scale).    We note that if the logistic map produced binary output strings with a uniform probability over strings, in theory the frequency would grow exponentially with complexity $r$, 
however, for Lempel-Ziv complexity with short strings, the distribution is not quite exponential, see \cite{dingle2018input,mingard2023deep}. 

The appearance of a distribution that is much less peaked towards high complexity than that of randomly chosen strings can be rationalised from AIT arguments.  To first order there are $\sim2^r$ strings of complexity $r$, each with probability $\sim2^{-r}$, such that the product  $2^r2^{-r}=O(1)$ is independent of $r$, and presumably uniform. For prefix complexity, the number of strings with complexity $r$ is slightly less than $2^r$, and hence for a prefix UTM the distribution would be biased to lower $r$ values to some extent.  Nonetheless, the brief rough argument outlined is still valid as a first order approximation for distributions with large spread, as we see in  Figure \ref{fig:distn_K}. See also ref.\ \cite{johnston2022symmetry} for similar arguments and results regarding the distribution of complexities in a biological setting, as well as ref.\  \cite{mingard2023deep} for the setting of machine learning. 

Finally, we note that in practice, far fewer samples are needed to produce a $P(\tilde{K}(x)=r)$ distribution than a $P(x)$ distribution, because many strings have the same $\tilde{K}$. Comparing the sampled distribution to a null-model of random strings may be the quickest and easiest way to diagnose simplicity bias in a dynamical system.

\subsubsection{Complex and pseudo-random outputs}
The occurrence of complex patterns may prompt a question: since the logistic map in Eq.\ (\ref{eq:logisticmap}) is itself simple with low Kolmogorov complexity, it might be supposed that any output $x$ it generates should be simple as well. So, how can we obtain large variations in complexity which would lead to large variations in probability? Or are all the apparently complex outputs merely pseudo-random patterns? We can address these analytically via bounding the complexity $K(x)$ of a discretised trajectory written as an $n$-bit binary string $x$ using the following inequality
\begin{equation}
K(x) \leq K(\text{logistic map}) + K(\mu,x_0) + K(n)+O(1)\label{eq:kx_bound}
\end{equation}
This bound follows from the fact that any $x$ can be described precisely by first describing: (a) the logistic map function in Eq.\ (\ref{eq:logisticmap}) with only a few bits because $K(\text{logistic map})$=$O(1)$ bits, (b) a pair of values  $\mu$ and $x_0$ which yield $x$ with $K(\mu,x_0)$ bits, and (c) the length of the string (i.e., the number of iterations to perform) via $n$, which is at most $\log_2(n)$ bits (up to loglog terms). From this upper bound, we see that if $\mu$ and $x_0$ are simple values with short descriptions, then $K(x)\ll n$ so that $x$ must be simple. However, because $\mu$ and $x_0$ are (truncated) real numbers randomly chosen from their intervals, they need not be simple --- rather the opposite --- because almost all decimal numbers sampled here are given to $d$ decimal places, and so will have high complexity values of $\sim\log_2(10)d$ bits. Therefore, outputs need not be simple just because the map is simple, due to the presence of the complexity of the inputs in Eq.\ (\ref{eq:kx_bound}). Nor are all outputs merely pseudo-random patterns. Indeed, this argument accords with the fact that in Figure \ref{fig:bias_with_simplicity bias} we see that very many outputs $x$ are realised, most of which must be (truly) complex because it is well known from AIT that almost all strings of length $n$ are complex and only relatively few are simple.

Extending the preceding discussion on Eq.\ (\ref{eq:kx_bound}), it is also interesting to ask if there are simple input parameters which generate high entropy pseudo-random outputs. Indeed there are, for example if  $\mu=4.0$ then almost all $x_0$ lead to chaotic trajectories including some simple initial values like $x_0=1/3$ which can be described with only a few bits, so that $K(\mu,x_0)=O(1)$ and so $K(x)\leq \log_2(n)+O(1)$ (up to additive loglog terms), but at the same time $\tilde{K}(x)\approx n$ bits because the Lempel-Ziv complexity measure cannot compress pseudo-random strings. 

On reflection, these pseudo-random strings with $K(x)\ll \tilde{K}(x)$ must be `rare' in the space, and not have individually high probability, otherwise we would not see simplicity bias, or at least see many points strongly violating of the upper bound in Figure \ref{fig:bias_with_simplicity bias}(a), which we do not. 

\subsubsection{Pre-chaotic regime}
As discussed above, when $\mu>3.5699...$ some trajectories can be chaotic, but there are also intermittent `islands' of stability, and in general the pattern of types of behaviour in the region 3.57 to 4.0 is intricate and complicated. However, in the region $\mu\in [0.0, 3.5699]$ the types of behaviour are more straightforward, with progressive period-doubling bifurcations but no true chaos.  Feigenbaum \cite{feigenbaum1979universal,feigenbaum1983universal} famously showed that the distance between successive bifurcations eventually decreases exponentially, as $\mu_{q+1}-\mu_{q}\approx 2.069/4.669^{q}$ \cite{berger2001chaos}, where $\mu_{q}$ is the value of the $q^{th}$ bifurcation.  

\begin{figure*}[htp]
\begin{center}
\subfigure[]{\label{fig:edge-a}\includegraphics[height=5.5cm,width=5.5cm]{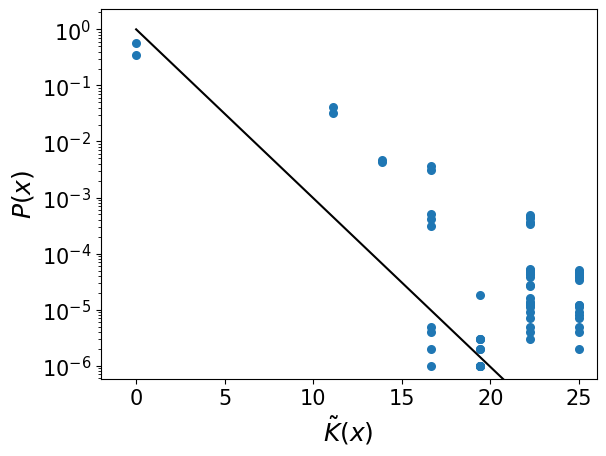}}
\subfigure[]{\label{fig:edge-a}\includegraphics[height=5.5cm,width=5.5cm]{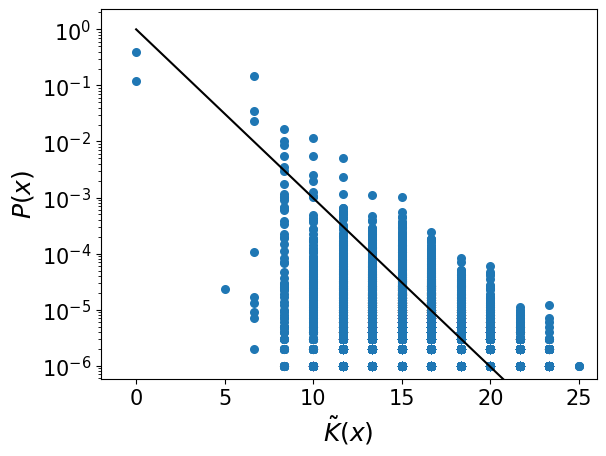}}
\subfigure[]{\label{fig:edge-a}\includegraphics[height=5.5cm,width=5.5cm]{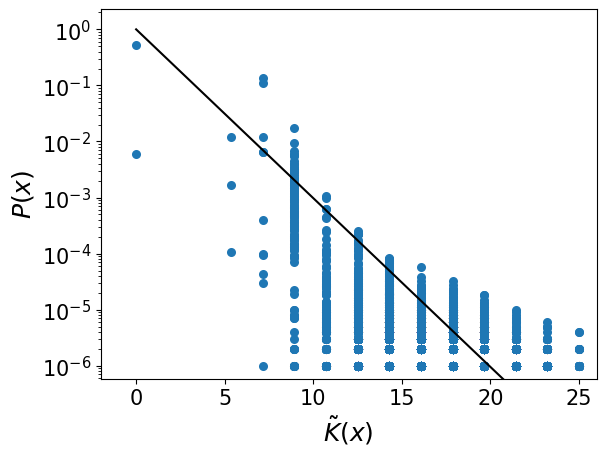}}
\end{center}
\caption{Simplicity bias in (a) the logistic map with $\mu$ sampled in [0.0, 3.5699], which is the non-chaotic period doubling regime (upper bound fitted slope is -0.17); (b) the Gauss map (upper bound fitted slope is -0.13); and (c) the sine map (upper bound fitted slope is -0.17). }
\label{fig:gauss_sine}
\end{figure*}

Because this interval $\mu\in [0.0, 3.5699]$ is relatively easy to understand, and to see if simplicity bias appears even without chaos, we also generated a complexity-probability plot for $\mu\in [0.0, 3.5699]$ uniformly sampled: Figure \ref{fig:gauss_sine}(a) shows the decay in probability with increasing complexity. As is apparent, there are  fewer patterns (blue dots) which is because we are restricting the possible space of patterns by restricting to this interval of $\mu$. We see that some simplicity bias is observed, but the bound is less clearly followed as compared to the when sampling across [0.0, 4.0]. For example, the upper bound on the data points' decay is less clearly linear.

To get a better understanding of the complexity of trajectories,  we can find\footnote{WolframAlpha (https://www.wolframalpha.com/), using command ``logistic map for r=[//number:3.4//] and x=[//number:0.2//]''} the oscillation periods for some chosen values of $\mu$: With $\mu=3.4$, period 2; $\mu=3.5$, period 4; $\mu=3.56$, period 8; $\mu=3.569$, period 32; $\mu=3.5699$, period 128. So for the sampled interval [0.0, 3.5699] the highest period is 128. Because we use $n=25$, any pattern with period 32 or more will not appear as periodic, because $n$ needs to be large for the periodicity to become plain. From this brief analysis, we can see how the low-probability and high-complexity patterns have appeared in Figure \ref{fig:gauss_sine}(a). In conclusion, some simplicity bias is observed for the  interval $\mu\in [0.0, 3.5699]$, but it is not as pronounced as for sampling $\mu\in [0.0, 4.0]$, which is presumably due to the presence of potentially more complex patterns in the interval 3.57 to 4.0.

\subsection{Gauss map (``mouse map'')}
 Moving on from the logistic map, we now explore simplicity bias in another 1D map, namely the \emph{Gauss map}, which is also known as the \emph{mouse map} because the bifurcation diagram is reminiscent of a mouse \cite{gaussmap}. 
The equation defining the dynamical system is 
\begin{equation}
x_{k+1} =e^{-\alpha x_{k}^2} + \beta \label{eq:gauss}
\end{equation}
The value $\alpha=7.5$ is chosen and will be fixed. For many values of $\alpha$ the bifurcation diagram is not sufficiently complex as to yield many varied trajectories; also the value 7.5 has been used previously  \cite{patidar2006co}. The value $\beta$ is the bifurcation parameter and the map is typically studied \cite{gaussmap,suryadi2020improvement} with $\beta\in[-1.0,1.0]$ . 

Similar to the logistic map example, we will sample initial values $x_0$ and values of $\beta$, then ignore the first 50 iterations (due to transient dynamics), and then digitise the real-valued trajectory to form binary string outputs. Because iterations of Eq.\ (\ref{eq:gauss}) are not confined to [0.0, 1.0] and due to the form of the bifurcation diagram \cite{gaussmap}, we will sample $x_0\in[-0.5, 0.5]$ uniformly. Also due to the form of the bifurcation diagram, the digitisation threshold will be set at 0.2 (instead of 0.5, as it was in  the case of the logistic map). Changing the threshold is done merely to avoid having too many trivial outputs of $x=000\dots000$. As above, we will use $n=25$.

In Figure \ref{fig:gauss_sine}(b), we see that there is clear simplicity bias in the Gauss map. The slope of the decay follows the predicted black line very closely, but the offset value $b$ is not well approximated by $b=0$. Nonetheless, there is clear simplicity bias, similar to the logistic map case. Hence, we can predict the relative log-probability of  output strings with quite high accuracy. 


\subsection{Sine map}
We now study simplicity bias in the \emph{sine map}. Wolfram \cite{wolfram_feigenbaum} described the sine map
\begin{equation}
x_{k+1} =  \mu\sin(\pi\sqrt{x_k})\label{eq:sine}
\end{equation}
and further displayed the bifurcation diagram for this map, which is broadly similar to that of the logistic map, to illustrate Feigenbaum's discovery of universality. 
For this map, we sample $x_0\in [0.0,1.0]$ uniformly $10^6$ times, and return to the digitisation threshold of 0.5. The parameter $\mu$ will sampled uniformly from [0.0, 1.0]. As before, the first 50 iterates will be ignored, and we will use $n=25$. Note that there is another form of the sine map which is sometimes studied in dynamical systems \cite{griffin2013sine,dong2021chaotification}, in which there is no square root on $x_k$.

Figure \ref{fig:gauss_sine}(c) shows the complexity-probability plot for the sine map, and simplicity bias is also present here. However, the upper bound is less clearly linear (on the log scale), and the black line prediction is not as closely followed. Also noteworthy is that in the tail of the decay at higher complexities, the upper bound  appears to curve up slightly. It is not clear why this happens.

\subsection{Bernoulli map}
Moving on to another prototypical chaotic map, the Bernoulli map (also known as the \emph{dyadic map}, \emph{bit shift map}, \emph{doubling map}, or \emph{sawtooth map}) is defined via the equation
\begin{equation}
x_{k+1} = ( 2x_k\mod 1)
\end{equation}
with $x_0\in (0.0, 1.0)$. This map shows sensitive dependence on initial conditions because the trajectories of two inputs $x_0$ and $y_0$ which differ by $|x_0 -y_{0}|\approx 0$  will eventually diverge for large enough $k$. 

Given that the Bernoulli map is also a 1D chaotic system with a simple ($O(1)$ complexity) map, it is interesting to ask whether it shows simplicity bias like the logistic map and others do. A little reflection shows that this map does not show simplicity bias nor bias, even for a digitised version. The trajectory is defined by multiplying 2 by $x_0$ ignoring the integer part, so if $x_0$ is a random real number then the trajectory will be random and incompressible because $(2x_0 \mod1)$ will be another random number, assuming that $x_0$ is defined to a large number of decimal places. Multiplying a random number by 2 does not remove the randomness. Hence the binary discretised trajectory sequence would look almost random, with small and quickly decaying autocorrelation (see also Section 9.4 of ref.\ \cite{calude2002information} for a similar conclusion).  For bias and simplicity bias, it is necessary for random inputs to be able to lose a lot of their information and complexity (i.e., complex inputs must be able to produce simple outputs), but the Bernoulli map does not allow this. Hence the behaviour of this map is similar to the logistic map with $\mu=4.0$ in the sense that there is no bias and no simplicity bias. Indeed, this similarity is quite natural due to the conjugacy between the logistic map with $\mu=4.0$ and the Bernoulli (doubling) map.
This map does not have a bifurcation parameter $\mu$.

\subsection{Tent map}
The last map we look at is the \emph{tent map}, which is quite well known and studied in dynamical systems research \cite{hasselblatt2003first}. The iterated values follow the function
 \[ 
 x_{k+1}  = \begin{cases} 
          2x_k &0 \leq x_k\leq 0.5 \\
          2-2x_k & 0.5< x_k\leq 1.0
       \end{cases}
 \]
with $x_k\in[0.0,1.0]$, for $k=0,1,2,3,\dots$. This map does not have a bifurcation parameter $\mu$. Despite being a 1D dynamical system, this map does will not lead to strong bias in the distribution of digitised binary string outputs $x$, and hence cannot possibly show simplicity bias. Intuitively, this can be seen due to the fact that almost all values of $x_0$ will yield complex paths,  while simplicity bias arises typically when most inputs lead to relatively simple paths. Indeed, because the tent map is topologically conjugate to the logistic map \cite{hasselblatt2003first} when $\mu=4.0$, and we saw neither bias nor simplicity bias in the logistic map \cite{hasselblatt2003first} when $\mu=4.0$, this helps to understand the absence of simplicity bias in the tent map.

\section{Discussion}
Arguments inspired by algorithmic information theory (AIT) predict that in many input-output maps, strongly non-uniform probability distributions over outputs will result, with complex patterns having exponentially low probabilities, and some simple patterns having high probability; this phenomenon is known as \emph{simplicity bias} \cite{dingle2018input}. Here, we numerically investigated the presence of simplicity bias in digitised trajectories arising from iterations of the logistic map, Gauss map, sine map, Bernoulli map, and tent map. By digitising the real-valued trajectories, we studied the probability and complexity of the resulting binary strings. Our main conclusions are that (i) we observe simplicity bias in the logistic map, Gauss map, and sine map, and also that in some cases the probability of resulting binary strings can be predicted \emph{a priori} with surprising accuracy; and (ii) we do not observe simplicity bias in the trajectories of the Bernoulli map and tent map, nor indeed any bias at all.  

Due to the qualitatively different behaviours exhibited by the logistic map for different $\mu$ values, we separately studied different regimes by sampling $\mu$ from (0.0, 4.0], (3.0, 4.0], (3.57, 4.0] and also $\mu=4.0$. In general,  simplicity bias and upper bound prediction accuracy was higher for $\mu$ sampled across the full range (0.0, 4.0] and  decreased for smaller ranges, until completely disappearing for $\mu=4.0$.  The logistic map is perhaps the most iconic example of a dynamical system in chaos theory, and has been very extensively studied for decades.  Here  we report a novel finding relating to this map, and one that is not (merely) a subtle higher-order effect, but rather a strong  effect related to order-of-magnitude variations in pattern probability.
This finding is also interesting given that we did not necessarily expect to observe simplicity bias when outputs can be pseudo-random (Cf.\ ref. \cite{dingle2018input}). It appears that in this map, pseudo-random outputs are sufficiently rare that they do not cause strong violations of the simplicity bias upper bound. Additionally, we found simplicity bias can be `tuned' via altering the  $\mu$ interval: sampling from the full interval (0.0, 4.0] yields a biased (low entropy) distribution over output strings along with simplicity bias, while sampling from higher values of $\mu\approx4.0$ yields low bias (high entropy) distributions and little or no simplicity bias.    
While we  observe simplicity bias similar to that predicted by AIT arguments in some of these maps,  we want to clarify that we do not think that we have in any way proven that these patterns are in fact directly linked to the AIT arguments.  For that, much more work is needed.   Nevertheless, we argue that studying probability-complexity relationships, and looking for patterns such as simplicity bias  may be a fruitful perspective for dynamical systems research.

A weakness of our probability predictions is that they only constitute an upper bound on the probabilities, and for example Figure \ref{fig:bias_with_simplicity bias}(a) shows that many output trajectory patterns $x$ fall far below their respective upper bounds. Following the hypothesis from \cite{dingle2020generic,alaskandarani2023low}, these low-complexity low-probability outputs are presumably patterns which the logistic map finds `hard' to make, yet are not intrinsically very complex. Further, the presence of these low-complexity low-probability patterns may indicate the non-universal power of the map \cite{dingle2020generic}. A potential avenue for future work would be studying what types of patterns occur far from the bound may, and possible approaches to improving on their probability predictions \cite{alaskandarani2023low}. 

The motivations for this work were to explore new examples of simplicity bias in physics and specifically dynamical systems, test the boundaries of relevance of simplicity bias, and explore how information theory and algorithmic probability can inform dynamical systems research. By extension, this work expands on the project of investigating the interaction between machine learning and dynamical systems because machine learning is intimately connected to information theory and algorithmics. 
The broader context of our work is a research project into the interface of dynamical systems, machine learning, and AIT. Given the strong interest in applying machine learning to the analysis of dynamical systems, an open and important question, then, is the extent to which machine learning methods are applicable in dynamical systems problems, and what kinds of limitations or advantages this relatively novel approach may have. Since information theory and machine learning are inextricably linked and have even been referred to as ``two sides of the same coin'' \cite{MacKay2003}, one perspective on these questions is to consider whether computation and information processing --- fundamental components of machine learning --- themselves might have limits of applicability, or in some other way constrain or inform dynamical behaviours. Some fascinating examples of such limits are known in relation to uncomputability of trajectories, meaning that some properties of dynamical trajectories may not be possible to calculate, even in principle. In a seminal paper, Moore \cite{moore1990unpredictability} proved that the motion of a single particle in a potential can (in some specific settings) be sufficiently complex as to simulate a Turing machine, and thereby yield trajectories with uncomputable quantitative properties. More recently, Watson et al. \cite{watson2022uncomputably} proved for an example many-body system that even if the exact initial  values of all system parameters were known, the renormalisaton group trajectory and resultant fixed point is impossible to predict. This kind of unpredictability is stronger than the unpredictability of (merely) chaotic dynamics, in which the limiting factor is the accuracy with which initial conditions can be measured.
See also Wolfram \cite{wolfram1985undecidability,wolfram2002new}, Svozil \cite{svozil1993randomness}, Lloyd \cite{lloyd2017uncomputability}, and Aguirre et al. \cite{aguirre2021undecidability} for more discussion of  (un)computability and (un)predictability in physical systems. Naturally, if the dynamics of some system cannot be computed even in principle, the accuracy of machine learning approaches to prediction in these settings will be restricted (but see \cite{lathrop1996learnability}).

There may also be deep connections between deep neural networks  (DNNs) and  simplicity bias. Indeed,  upon random sampling of parameters, DNNs exhibit an exponential bias towards functions with low complexity~\cite{valle2018deep,mingard2019neural,bhattamishra2022simplicity}. This property implies that they can learn Kolmogorov simple data fairly well, but will not generalise well on complex data.  Interestingly, by changing the initialization over parameters towards a more artefactual chaotic regime of DNNs, this simplicity bias becomes weaker~\cite{yang2019fine}.  It has been recently shown that in this regime,  DNNs no longer generalise well on both simple and complex data~\cite{mingard2023deep}, and tend to overfit.  This is not unexpected, because  DNNs are highly expressive, and classical bias-variance arguments suggest that they should be highly prone to overfitting.    The big question is why standard DNNs don't fall prey to this problem.   It was argued  that the  Occam's razor like  \emph{simplicity bias}  toward simple functions observed in standard DNNs compensates the exponential growth of the number of possible functions with complexity, and that this compensation explains why such DNNs can be highly expressive without overfitting \cite{mingard2023deep}.  These principles imply that if  a dynamical system exhibits some form of simplicity bias then the inbuilt Occam's razor inductive bias of DNNs should make it much easier to learn by DNNs than in the opposite case where the dynamical system does not have simplicity bias. 


The word ``complexity'' can take on many meanings \cite{lloyd2001measures,mitchell2009complexity}, and can be vague. In this work we are precise about what we mean by complexity, which is Kolmogorov complexity, and in practice using lossless compression methods which is a standard and theoretically motivated approximation to the true uncomputable quantity. Bialek et al. \cite{bialek2001complexity,bialek2001predictability} discuss complexity in relation to time series, and argue that Kolmogorov complexity lacks in its intuitive appeal for a measure of the complexity of these types of sequential patterns. Many would agree that a series with rich complicated structure and long-range correlations is truly complex, whereas a random string of bits is merely an irregular and in a sense trivial pattern. In contrast, Kolmogorov complexity assigns the highest complexity to such random strings, precisely because they do not contain any correlations or structure. Having noted this, the discussion does not directly bear upon our work, because we are not studying `complexity' in a general sense or trying to argue for one or other metric. Rather, we are studying specifically simplicity bias and AIT inspired  arguments as a mathematical framework for understanding probability and dynamical systems. AIT allows one to make quantitative bounds and predictions about various systems, as we illustrate here, regardless of whether or not the term ``complexity'' is being used in it truest or most correct sense.

Mathematicians and physicists are fascinated by studying simple systems which can exhibit complex patterns, even while they are not technically random. In this context, Wolfram \cite{wolfram2002new} investigated simple automata, and showed that some have high levels of computing power, leading to his conjecture that many natural systems can be Turing complete, i.e., that they can implement arbitrary algorithms, and hence produce complex patterns. Despite this focus on complexity, we observe in our work here that complexity is in a sense not actually that common: Even though the logistic map is famous precisely due to its ability to produce  complex behaviour, within the space of possible parameters, i.e., $\mu,x_0\in\mathcal{R}$ and even restricting to $\mu\in(0.0,4.0]$ and $x_0\in(0.0,1.0)$, chaos and `complexity' only occur rarely (assuming uniform sampling of the parameter space). This observation accords with the fact that while Wolfram highlighted some rule sets that produce randomness, most of his automata rules sets do not produce complex pseudo-random patterns \cite{wolfram2002new}. Coe et al.\  \cite{coe2008cellular} analytically studied the question of when automata produce `random' and complex dynamics, and also found that relatively few are random.

The study of random dynamical systems --- in which dynamical systems are randomly perturbed in some way --- is quite well established \cite{arnold1995random}, including specifically dynamical systems with additive noise \cite{doan2018hopf}. While deterministic dynamical systems have been studied due to their relevance to modelling problems in e.g., ecology, physics, and economics, randomly perturbed dynamics are also common in science and hence important to study, and arise naturally in modelling physical systems (e.g., \cite{dingle2012knudsen}). This motivates studying simplicity bias also in such random dynamical systems in future work. As a separate motivation, in the deterministic logistic map, many sampled $\mu$ values yield trajectories which quickly converge to fixed points or other fairly trivial patterns. However, by introducing some small random noise these trajectories may be prevented from converging into trivial patterns, and therefore may show simplicity bias even while the deterministic counterpart does not. The relation of simplicity bias to random dynamical systems has been initially studied recently \cite{hamzi2024simplicity}, but many questions remain open.

In this work we have used 1D maps including the logistic map as toy systems to explore simplicity bias in dynamical systems. The connection between AIT and dynamical systems has  received some earlier attention from an analytical perspective \cite{white1993algorithmic,brudno1978complexity,v2022ergodic}, and computational perspective \cite{zenil2019algorithmic,hamzi2024simplicity}. In future work it may be fruitful to investigate which, if any, properties of the logistic (or other chaotic) map can be predicted or bounded using the simplicity bias bound. That is, if we assume that the simplicity bias bound holds, what other dynamical properties of the map might follow from this? Even if these properties are already known, it would still be insightful to see if they could be derived (perhaps more simply) from AIT arguments. Further, the presence of simplicity bias in some of the 1D map may lead to using the simplicity bias bound as  a trajectory prediction method. This is related to the \emph{a priori} prediction of natural time series patterns attempted by Dingle et al.\ \cite{dingle2023note}. Another angle would be to study non-digitised trajectories of dynamical systems, which would require different complexities measures amenable to continuous curves, such as proposed in ref. \cite{terry2023fourier}.
More generally, exploring the use of  arguments inspired by AIT in dynamical systems and chaos research is an attractive and potentially fruitful research avenue.

\vspace{0.5cm}
\noindent
{\bf Acknowledgments:} This project has been partially supported by Gulf University for Science and Technology, including by project code: ISG Case 9. We  thank G. Valle-Perez for valuable discussions.

\bibliographystyle{unsrt}
\bibliography{SB_1Dmaps_Refs} 


\newpage

\appendix

\section{The impact of the number of iterations}\label{numbiterations}

In this brief section, we give examples of the impact of changing the length $n$ of binary strings.

\begin{figure*}[htp]
\begin{center}
\subfigure[]{\label{fig:edge-a}\includegraphics[height=5.5cm,width=5.5cm]{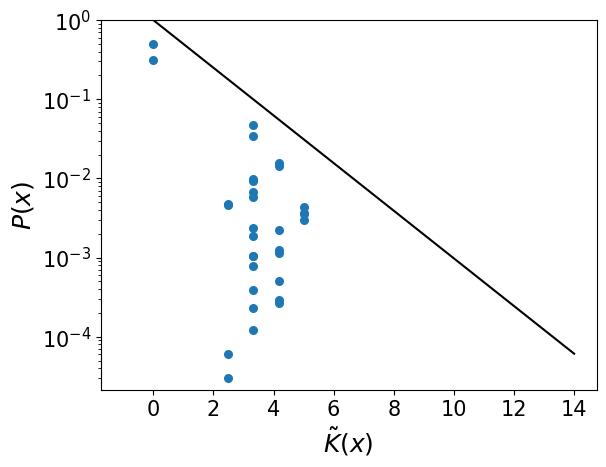}}
\subfigure[]{\label{fig:edge-a}\includegraphics[height=5.5cm,width=5.5cm]{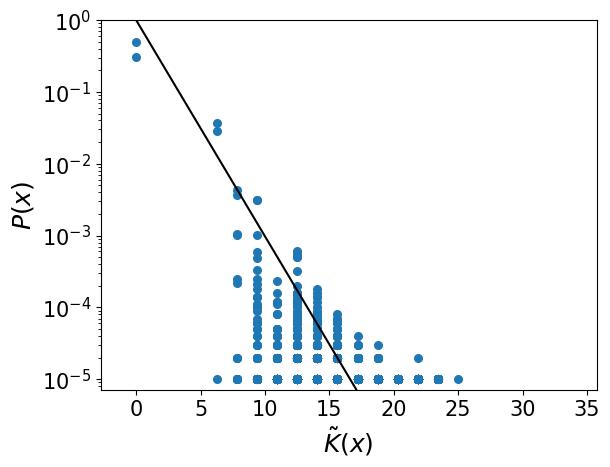}}\\
\subfigure[]{\label{fig:edge-a}\includegraphics[height=5.5cm,width=5.5cm]{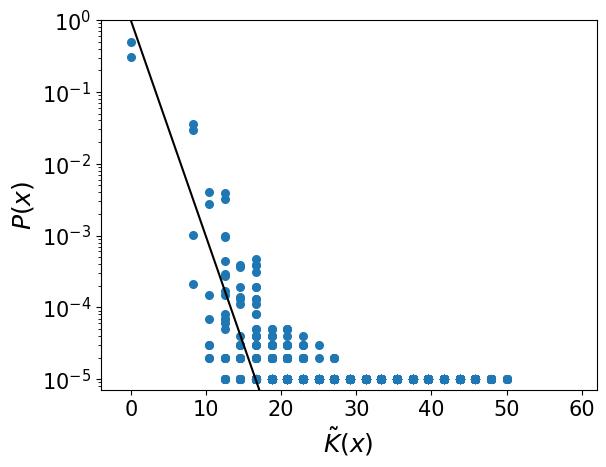}}
\subfigure[]{\label{fig:edge-a}\includegraphics[height=5.5cm,width=5.5cm]{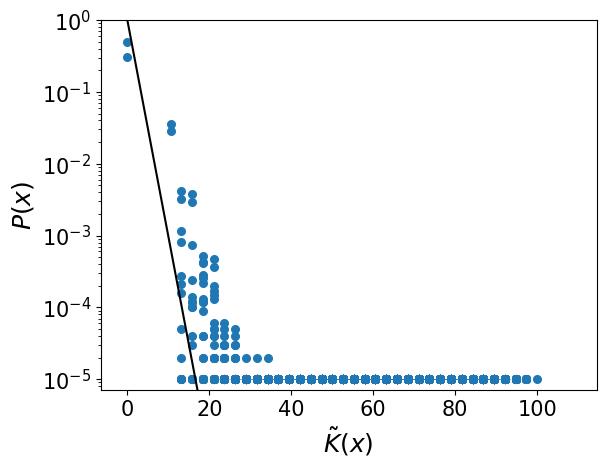}}
\end{center}
\caption{Simplicity bias with different number of iterations. (a) With $n=5$ iterations, there is some simplicity bias but it is not pronounced; (b) with $n=25$ iterations, simplicity bias is very clear; with (c) $n=50$ iterations there is still clear simplicity bias, but a long `tail' begins to emerge, of low frequency patterns; (d) with $n=100$ iterations, there is still some simplicity bias but the `tail' has become more dominant and simplicity bias is less clear.  }
\label{fig:niterations}
\end{figure*}

For very small numbers of iterations, simplicity bias is not pronounced. This is due primarily to the fact that the complexity measure is not very sensitive to small changes in complexity values. To illustrate, in Figure \ref{fig:niterations}(a) a plot is shown for $n=5$ iterations only. In this case, there is some evidence of simplicity bias (i.e., an inverse relation between complexity and probability), but the relation is not strong. A similar weak connection between complexity and probability was  observed for very short sequences in an earlier time series study of simplicity bias \cite{dingle2023note}. In Figure \ref{fig:niterations}(b), clear simplicity bias is observed with $n=25$. At this length, the complexity measure is sufficiently sensitive and we have sufficient sampling to reveal order of magnitude variations in probability. 

In theory, for larger output complexities simplicity bias should be more clear \cite{dingle2018input}, and in fact AIT results in general are typically stated up to $O(1)$ terms, which become irrelevant  when pattern complexities become asymptotically large \cite{li2008introduction}. However, in this logistic map example we find that for longer iterations, the pattern of simplicity bias becomes more noisy and less clear; see Figures \ref{fig:niterations}(c) and (d). For larger $n$ the number of possible patterns becomes exponentially large and the one million samples produces many patterns of frequency $\approx$1, yielding a long flat `tail' in the graph.

\section{Alternate figures highlighting density of points}
In this section we repeat the some of the plots from the main text, but now giving semi-transparent datapoints, so as to highlight the density of data points.

\begin{figure*}[htp]
\begin{center}
\subfigure[]{\label{fig:edge-a}\includegraphics[height=5.5cm,width=5.5cm]{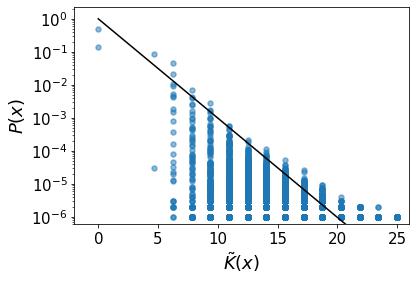}}
\subfigure[]{\label{fig:edge-a}\includegraphics[height=5.5cm,width=5.5cm]{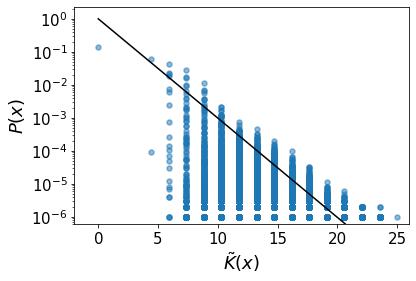}}\\
\subfigure[]{\label{fig:edge-a}\includegraphics[height=5.5cm,width=5.5cm]{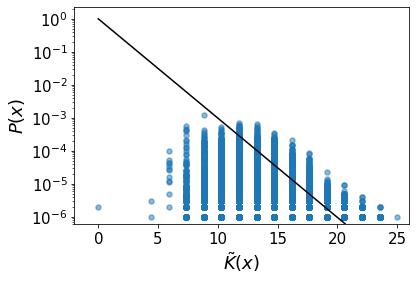}}
\subfigure[]{\label{fig:edge-a}\includegraphics[height=5.5cm,width=5.5cm]{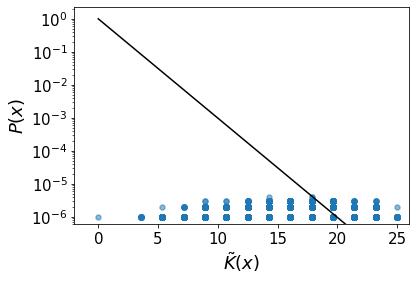}}
\end{center}
\caption{Simplicity bias in the logistic map, same as Figure 3, but with semi-transparent data points.}
\label{fig:app_transparent_3}
\end{figure*}

\begin{figure*}[htp]
\begin{center}
\subfigure[]{\label{fig:edge-a}\includegraphics[height=5.5cm,width=5.5cm]{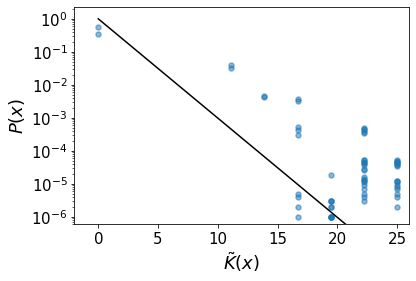}}
\subfigure[]{\label{fig:edge-a}\includegraphics[height=5.5cm,width=5.5cm]{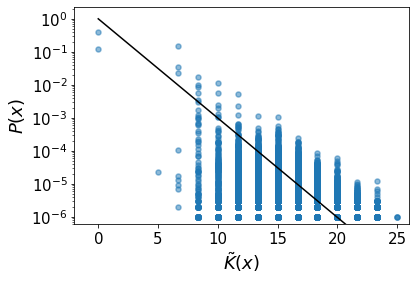}}\\
\subfigure[]{\label{fig:edge-a}\includegraphics[height=5.5cm,width=5.5cm]{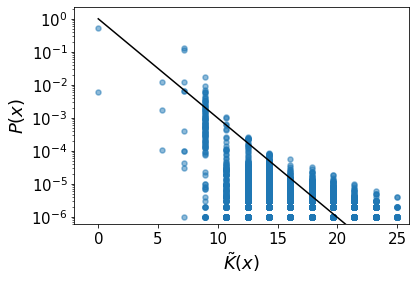}}
\end{center}
\caption{Simplicity bias in the logistic, Gauss map, and sine map, same as Figure 5, but with semi-transparent data points.}
\label{fig:app_transparent_5}
\end{figure*}

\end{document}